\documentclass{article}
\usepackage[margin=25mm]{geometry}
\usepackage{geometry}
 \usepackage[english]{babel}
\usepackage[square,sort,comma,numbers]{natbib}
\usepackage{xcolor}
\usepackage{lineno,hyperref}
\modulolinenumbers[5]
\usepackage{amsmath, amssymb, tabularx, amsfonts,dirtytalk}
\usepackage{tikz}
\usetikzlibrary{trees}
\usepackage{graphicx} 
\usepackage{url,microtype}
\usepackage{longtable}
\usepackage{wrapfig}
\usepackage{verbatim}
\usepackage[labelfont=bf]{caption}
\usepackage{multirow}
\usepackage{soul}
\usepackage{color}
\usepackage{soul}
\usepackage{color}
 \usepackage{xspace}
\usepackage{xparse}
\NewDocumentCommand{\cbox}{s O{1ex}}{%
  \setlength{\fboxsep}{-\fboxrule}%
  \IfBooleanTF{#1}
    {\frame{\rule{0.5\dimexpr#2}{0.5\dimexpr#2}\rule[0.5\dimexpr#2]{0.5\dimexpr#2}{0.5\dimexpr#2}}}
    {\fbox{\rule[0.5\dimexpr#2]{0.5\dimexpr#2}{0.5\dimexpr#2}\rule{0.5\dimexpr#2}{0.5\dimexpr#2}}}
  \xspace
}
 \usepackage{comment}

 \usepackage{soul,color}
\usepackage{multirow}

\usepackage{multirow}
\usepackage{caption}

\usepackage[utf8]{inputenc}
\usepackage{longtable}
\usepackage{amsmath, amssymb, tabularx, amsfonts,dirtytalk}

\usepackage{booktabs}
\usepackage{caption}
\usepackage{siunitx}
\usepackage{float}

\title{Scheduling a Multi-Product Pipeline:  
\\ A Discretized MILP Formulation}
\author{Ales Wodecki\thanks{These authors contributed equally to this work.}, Pavel Rytíř\footnotemark[1], Vyacheslav Kungurtsev, Jakub Marecek\footnote{All authors are at the Czech Technical University.}}
\date{September 2023}

\begin{document}

\maketitle

\begin{abstract}
Multi-product pipelines are a highly efficient means of transporting liquids. Traditionally used to transport petroleum, its products and derivatives, they are now 
being repurposed to transport liquified natural gas admixed with hydrogen of various colors. We propose a novel mixed-integer linear programming (MILP) formulation, which optimizes efficiency while satisfying a wide range of real-world constraints. The proposed formulation has been developed to meet the needs of the Czech national pipeline operator \v{C}EPRO, who uses it in monthly planning of operations. We provide tests on well-known synthetic (path-graph) networks and demonstrate the formulation's scaling properties using open-source and commercial MILP solvers.
\end{abstract}

\section{Introduction}
The interest in logistics problems associated with pipeline transport is due to the inherent efficiency of this mode of transportation \cite{Li2021, arya2022recent},
and wide-spread availability of pipeline transport systems. 
Multi-product pipelines are of particular interest, since they provide additional flexibility compared to single commodity pipelines \cite{arya2022recent}. Furthermore, most existing networks are of the multi-product variety \cite{LI2021395, TU202387, Li2021}. These include conventional pipelines for the transport of petroleum and its derivatives, as well as pipelines adapted to transport of liquified natural gas (LNG) and hydrogen \cite{Li2021, TIAN2023e21454, Zhang23}. 
(In transporting hydrogen, one may wish to distinguish green hydrogen, which is made by electrolysers powered by renewable energy sources, from the rest of the hydrogen color spectrum.
Before pure hydrogen is transported, various mixtures of LNG and hydrogen will be transported.)
To solve the problem of scheduling such a pipeline, many different solution strategies have been applied in the past including mixed integer linear programs (MILP), heuristics, and neural networks. Solutions in which the aforementioned methods are combined have also been studied \cite{XU2021341, BAMOUMEN2023106082}. In general, neural networks and heuristics provide fast solutions for some instances, but remain a black box, without satisfactory guarantees in general \cite{arya2022recent, YUCEER2005933, BAMOUMEN2023106082, Liang16}. The use of MILP solvers guarantees convergence at the cost of a possibly longer runtime.

In the literature, MILP models are generally classified according to pipeline topology, optimization objective, and means of uncertainty quantification \cite{Li2021}. 
Here, we focus on a demand-driven scenario on a path graph and assume that there are point estimates for all quantities. 
There is a long history of work on this scenario. 
Recently, 
Rejowski et al.\ \cite{REJOWSKI20081042} have provided a pure MILP formulation for this scenario and have demonstrated adequate scaling properties on a 5 node network with a time horizon of 3 days. 
Meira et al.\ \cite{meira20} have combined MILP techniques with a heuristic that uses a fitness function to evaluate the quality of solutions. 
This combination can solve instances involving up to 6 sites with a relatively long time horizon of 30 days. 
In parallel, Xu et al.\ \cite{XU2021341} proposed a decomposition method, which allows solving scheduling problems at up to 6 sites over a 14-day time horizon. Note that all of the aforementioned solutions were obtained within an hour per instance.

We provide a novel MILP formulation of the multi-product demand-driven pipeline scheduling problem, in which the scheduling of product transport is treated as a bin-packing problem \cite{martello1990knapsack}, using discretization of both time and space \cite{Allen_2012}. 
Such formulations are sometimes known as discretized \cite[e.g.]{CORREIA20082103} or space-indexed \cite[e.g.]{Allen_2012}. 
A bevy of operational constraints is added to the bin-packing formulation to arrive at the model of a multi-product pipeline. 
A branch-and-cut algorithm is developed in SCIP, which provides a speedup for certain problem instances. 
We also show that the formulation scales well using Gurobi Optimizer, the state-of-the-art commercial solver, without the need for custom separation of cutting planes. 
Using SCIP with the branch-and-cut algorithm, we can scale up to 8 sites and 576 time units (24 days) within 1000 seconds of run-time on an Apple Macbook Pro.
Using Gurobi Optimizer, we can scale up to 12 sites and 744 time units (31 days), with optimality gap below 10$^{-4}$ after 1 hour of runtime on the same Apple Macbook Pro. 
Finally, we note that the model has been developed in cooperation with the Czech national pipeline operator \v{C}EPRO, who currently utilizes the formulation to provide monthly plans for the operations of their multi-product pipeline system.

\section{MILP Formulation}\label{sec_MILP_form}

This section presents the full MILP formulation. When referring to the graph representing the multi-product pipeline system, we often call its vertices sites, and its edges pipes or pipe segments. The sites may be storage sites or refineries. We keep track of the capacity only at the storage sites, and the transport between sites is realized using pumping regimes that in general may use several pipes of the network and numerous pumps. Each pumping regime is identified by its source, destination, and the pipes utilized. 
The primary objective of the operator is to satisfy demand at the storage sites, while maximizing the inflow of product into the network from the refineries, up to a limit given by refineries' operating constraints. (The refineries' operating constraints remain private to the refineries.) Further terms may also be added to the objective, such as minimizing pumping costs or accommodating preference regarding the distribution of the products across the sites at the end of the planning horizon (see Section \ref{sec_structure_of_obj}). First, we present a short overview of the notation we use in Table~\ref{tab:notation}.

\begin{longtable}[tb]{lcp{14cm}} 
\caption{Table of notation}\\    
\label{tab:notation}
$\mathcal{S}$ & :& The set of all sites in the network, this includes storage sites and refineries\\
$\mathcal{S}_{\text{store}}$ & :& The set of all storage sites $\mathcal{S}_{\text{store}} \subset \mathcal{S}$\\
$\mathcal{S}_{\text{ref}}$ & :& The set of all refineries $\mathcal{S}_{\text{ref}} \subset \mathcal{S}$\\
$\mathcal{P}$ & : & The set of products\\
$\mathcal{P}_{\text{flush}}$ & : & The set of all flushing products $\mathcal{P}_{\text{flush}} \subset \mathcal{P}$\\
$\mathcal{P}_{\text{stain}}$ & : & The set of all staining products $\mathcal{P}_{\text{stain}} \subset \mathcal{P}$\\
$\mathcal{T}$ & : & Discrete time horizon over which the problem is solved\\
$\mathcal{R}$ & : & The set of pumping regimes
\\
$\mathcal{E}$ & : & The set of edges (pipelines), each of which is represented by a pseudo bin packing problem\\
$\mathcal{B}_{s, p}^{V}$ & : & Volumes of batches containing product $p$ that originate at site $s$
\\
$\mathcal{B}_{s,p,r}^{V}$ & : & Volumes of batches containing product $p$ that originate at site $s$ associated with regime $r$
\\
$e_{\text{init}}$ & : & Denotes the set of all initial batches for $e \in \mathcal{E}$
\\
$e_{\text{transit}}$ & : & Denotes the set of all transit batches for $e \in \mathcal{E}$
\\
$e_{\text{final}}$ & : & Denotes the set of all final batches for $e \in \mathcal{E}$
\\
$r_{\text{orig}}$ & : & Denotes the first pipe which the pumping regime $r \in \mathcal{R}$ uses
\\
$r_{\text{dest}}$ & : & Denotes the final pipe which the pumping regime $r \in \mathcal{R}$ uses
\\
$\mathcal{F}_{r}$ & : & The flow rate of a given pumping regime $r \in \mathcal{R}$
\\
$\mathcal{E}_{r}$ & : & The ordered list of edges that the pumping regime $r \in \mathcal{R}$ utilizes
\\
$\text{pass}\left(e,r\right)$ & : & The amount of timestamps needed for a batch to pass through pipe $e$ when pumping regime $r$ is engaged
\\
$b_{p}$ & : & The product contained in batch $b$
\\
$b^{V}$ & : & Units of volume contained in batch $b$

\end{longtable}

\subsection{General Considerations}\label{General_section}
The transport through the pipeline network is optimized over a finite time horizon (typically taken to be a couple of weeks). This time horizon is discretized using a time step $\Delta t$ (a common choice is $\Delta t = 1 h$) and the resulting set of timestamps is denoted $\mathcal{T}$. Let $\mathcal{P}$ be the set of products and $\mathcal{B}_{s, p}^{V}$ be the possible batch sizes (in units of volume) of product $p \in \mathcal{P}$ originating from site $s \in \mathcal{S}$. The set $\mathcal{B}_{s, p}^{V}$ always contains a standard batch size for all $s$ and $p$.

In addition to the standard batch sizes, other batch sizes are introduced due to the following considerations. In multi-product networks, it is common to determine a default product that is left in the pipes when the pumps are idle. We call this product the flushing product, all other products are labeled as staining. The set $\mathcal{P}$ is a union of flushing products $\mathcal{P}_{\text{flush}}$ and staining products $\mathcal{P}_{\text{stain}}$. If $p\in\mathcal{P}_{\text{stain}}$ the batch size $\mathcal{B}_{s, p}^{V}$ is always a singleton that contains only the standard batch size. On the other hand if $p\in\mathcal{P}_{\text{flush}}$ the set $\mathcal{B}_{s, p}^{V}$ may contain more than one batch size. This is because the amount of product needed to fill the relevant pipe segments with the flushing product may be larger than the standard batch size. Pumping in the network begins by the activation of a pumping regime. These can possibly span several edges (pipes) of the network. A pumping regime that spans multiple edges originating at site $s$ can result in $\mathcal{B}_{s, p}^{V}$ being augmented by a large enough flushing batch that is able to clean the corresponding set of pipe segments. Details on the construction of the set $\mathcal{B}_{s, p}^{V}$ can be found in Section \ref{sec_flushing}.

These batch sizes are, together with the pumping regime parameters, used to determine the amount of transported mass units per batch as well as the time needed to complete the transport of a single batch. A pumping regime $r \in \mathcal{R}$ connecting sites $r_{\text{orig}}=s_{1},r_{\text{dest}}=s_{2}\in \mathcal{S}$ through edges $\mathcal{E}_{r} = \left(e_{1},\ldots,e_{k}\right)$ has a flow rate indicated by $\mathcal{F}_{r}$ (in volume units per $\Delta t$). Applying a regime $r$ such that $r_{\text{orig}}=s$ to a batch of a given volume $w\in\mathcal{B}_{s, p}^{V}$ determines the time it takes to pump the product into the pipeline (by $\mathcal{F}_{r} / w$). In the discretized time horizon $\mathcal{T}$, this procedure results in a batch length $L\left(r,p;j\left(w\right)\right)$, where $j\left(w\right)$ denotes the index of the batch. i.e. 
\begin{equation}\label{length_def}
L\left(r,p;j\left(w\right)\right)=\left\lceil \frac{\mathcal{F}_{r}}{w\left(r,p,r_{\text{orig}}\right)}\right\rceil.
\end{equation}
The next section is dedicated to showing how batch scheduling may be viewed as a bin-packing problem.

\subsection{Modeling the Transport of Product Through the Network as a Bin Packing Problem}\label{sec_bin_pack}
To model transport through the edges of the graph, we consider a one-dimensional bin-packing problem for each of the edges $e\in\mathcal{E}$ with the size equal to the cardinality of $\mathcal{T}$. Packing a box at time $t\in\mathcal{T}$ then corresponds to scheduling a batch at time $t$. To consider every possible scheduling decision for the bin packing corresponding to edge $e\in\mathcal{E}$, we use  $L\left(r,p;j\left(w\right)\right)$ from \eqref{length_def} to define the box length for each relevant tuple $\left(r,p;j\right)$. A tuple is relevant for a bin packing problem $e\in\mathcal{E}$, if the first element in the ordered tuple $\mathcal{E}_r$ is $e$. This corresponds to the situation in which pumping is initiated at edge $e$. A tuple may also be relevant if the batch enters the pipe $e$ as a result of a pumping regime initiated at another edge; this is explained in detail in Section \ref{sec_pregs_routes}. To denote the set of relevant tuples, the following notation is used, when there is no risk of confusion:
\begin{equation}
    \left(r,p;j\right)=b\in e\in\mathcal{E},
\end{equation}
where the edges $e$ are identified with the bin packing and hold all of the relevant tuples, which are labeled $b$ for simplicity. To adapt the terminology to current use, we call the packed boxes $b$ batches.

A discretized formulation of the bin-packing problem \cite{Allen_2012} considers a binary decision variable, which is 1 when a batch is scheduled, which corresponds to a box being placed in terms of the bin packing. More succinctly,

\begin{equation}
            v_{ebt}  
        \begin{cases}
        1, & \text{if a batch corresponding to b = $\left(r,p;j\right)$ starts at time $t\in\mathcal{T}$ in edge $e \in \mathcal{E}$}, \\
        0, & \text{otherwise.}
        \end{cases}
\end{equation}
The relevant constraints read
 
\begin{equation}
    \sum_{b\in e}v_{ebt}\leq1, \indent \forall e\in\mathcal{E}, \forall t \in \mathcal{T}, \label{box_spatial_count_1}
\end{equation}

\begin{equation} \label{eq3}
\begin{split}
t +  L\left(b\right) v_{ebt} \leq \mathcal{T}_{\text{max}}, \indent \forall e\in \mathcal{E}, \forall t\in \mathcal{T}, \forall b = \left(r,p;j\right),
\end{split}
\end{equation}

\begin{equation}
    \sum_{\left(b^{\prime},t^{\prime}\right)\in\text{Excl}\left(e,t\right)}v_{eb^{\prime}t^{\prime}}\leq1,  \indent \forall e\in \mathcal{E}, \forall t\in \mathcal{T}, \label{box_overlap_exclusion}
\end{equation}

where $\text{Excl}\left(e,t\right)$ is a set that contains pairs of decision variables that must be mutually excluded to prevent the overlap of batches (see \cite{Allen_2012} for details). Notice that compared to \cite{Allen_2012} the constraints for proper batch scheduling are a bin-packing problem with the upper bound for batch placement removed. Thus, we refer to (\ref{box_spatial_count_1})-(\ref{box_overlap_exclusion}) as the pseudo-box packing problem.

\subsection{Pumping Regimes and Routes\label{sec_pregs_routes}}
In the previous section, the scheduling of batches using a bin-packing problem was discussed. Such scheduling, however, is not a matter of a single edge and thus a single pseudo-bin packing problem. Every transport through the network is a consequence of the activation of a pumping regime $r \in \mathcal{R}$ that spans the edges $\mathcal{E}_{r} = \left(e_{1},\ldots,e_{k}\right)$. As discussed in Section \ref{General_section}, a pumping regime $r \in \mathcal{R}$ spanning the edges $\mathcal{E}_{r} = \left(e_{1},\ldots,e_{k}\right)$ may be applied to a product $p \in \mathcal{P}$ to give rise to a batch $b = \left(r,p;j\right)$ that is to be placed in the bin packing problem of edge $e_{1}$. If the cardinality of $\mathcal{E}_{r}$ is greater than 1 then the placement of a batch in the bin packing problem $e_{1}$ must result in the placement of a batch (with the same length) in problem $e_{2}$ which expresses the movement of the batch through the pipeline. 

More precisely, a batch $b$ in a bin-packing problem corresponding to edge $e$ always belongs to a pumping regime i.e. $b = \left(r,p;j\right)$ and $e \in \mathcal{E}_{r}$ for some $r\in\mathcal{R}$. Let $b\left(e_{1};p,r;j\right),\ldots,b\left(e_{k},p,r;j\right)$ be batches corresponding to the realizations of a regime $r$ applied to product $p$. The route constraints can then be formulated as follows

\begin{equation}
    v_{e_{i}b\left(e_{i};p,r\right)t}=v_{e_{i+1}b\left(e_{i+1};p,r\right)t} \quad \forall i\in\left\{ 1,\ldots,k-1\right\} ,t\in\mathcal{T}.
    \label{routes_constrains}
\end{equation}
In the following sections, some additional constraints relating to batch transport and tank capacity modeling are explained. The constraints (\ref{routes_constrains}) imply that the behavior of the whole chain of decision variables is determined by the placement of batch on the first edge of the regime. Thus, it is useful to adapt the following terminology. We call $b\left(e_{1};p,r\right)$ an initial batch, $b\left(e_{i};p,r\right)$ for $i \in \left\{ 1,\ldots,k-1\right\} $ a transit batch and $b\left(e_{k};p,r\right)$ a final batch. The sets of all initial, transit and final batches of an edge $e$ are denoted $e_{\text{init}},e_{\text{transit}}$ and $e_{\text{final}}$ respectively. This division comes in handy when formulating tank constraints in the following sections.

\subsection{Flushing and Staining\label{sec_flushing}}
As mentioned in Section \ref{General_section}, the pipeline must be filled with a predetermined product when idle. This product is labeled flushing, and any other product is labeled staining. When a route is idle, it is filled with a flushing product. In the following, the constraints that enforce this behavior are described. 

Let $\mathcal{P}_{\text{stain}} \subset \mathcal{P}$ label the subset of products that stain the pipeline and let $\mathcal{P}_{\text{flush}} \subset \mathcal{P}$ be the subset of products that can flush the pipeline. All the products considered fall into one of these categories, i.e.,

\begin{equation}
    \mathcal{P}_{\text{stain}} \cup \mathcal{P}_{\text{flush}} = \mathcal{P}.
\end{equation}
The flushing mechanic excludes illegal sequences of products (two different staining products after one another) and makes sure that a pipeline can only be left idle while containing the flushing product. Based on the volume of the pipeline that needs to be filled with the flushing product an additional batch size needs to be added in order to ensure that the relevant pipe segments are filled with flushing product. This volume is regime-specific and denoted $r_{V}$.

Once this volume is determined, the set of all possible batches for product $p$ originating at a site $s$ $\mathcal{B}_{s, p}^{V}$ is defined as

\begin{equation}
\begin{array}{cc}
\mathcal{B}_{s,p}^{V}=\left\{ B_{s,p}^{\text{standard}}\right\}  & \text{\ensuremath{\quad\text{for }p\in\mathcal{P}_{\text{stain}}}},\\
\mathcal{B}_{s,p}^{V}=\left\{ B_{s,p}^{\text{standard}}\right\} \cup\left\{ r_{\text{V}}:r\in\mathcal{R},r_{\text{orig}}=s,r_{\text{V}}>B_{s,p}^{\text{standard}}\right\}  & \text{\ensuremath{\quad\text{for }p\in\mathcal{P}_{\text{flush}}}},
\end{array}\label{batches_def}
\end{equation}
where $B_{s,p}^{\text{standard}}$ is the standard batch size defined for each site $s$ and product $p$. Additionally, define

\begin{equation}
    \begin{array}{cc}
\mathcal{B}_{s,p,r}^{V}=\left\{ B_{s,p}^{\text{standard}}\right\}  & \text{\ensuremath{\quad\text{for }p\in\mathcal{P}_{\text{stain}}}},\\
\mathcal{B}_{s,p,r}^{V}=\left\{ B_{s,p}^{\text{standard}}\right\} \cup\left\{ r_{\text{V}}:r_{\text{V}}>B_{s,p}^{\text{standard}}\right\}  & \text{\ensuremath{\quad\text{for }p\in\mathcal{P}_{\text{flush}}}},
\end{array}
\end{equation}
as the regime specific analogue of (\ref{batches_def}) i.e. $\mathcal{B}_{s,p,r}^{V}$ contains only the batches associated with regime $r$.

Next, the flushing rules are outlined from a birds-eye view, after which the succinct formulation follows. Flushing is regime-specific; this means that whenever a regime pumps a product $p\in\mathcal{P}_{\text{stain}}$, the same regime needs to be used to flush. This ensures that the entire route is left in a nonstained state after a sequence of operations involving one pumping regime is performed. Whenever a product $p\in\mathcal{P}_{\text{stain}}$ is pumped using regime $r$, only two followup actions are permitted:
\begin{itemize}
    \item either the same product can be pumped again (using the same regime)
    \item or one of the flushing products can be pumped (using the same regime and sufficient batch size).
\end{itemize}

This mechanism is implemented in the MILP formulation by adding the following decision variables and constraints. First, the end points of batches that correspond to staining products are marked on each edge

\begin{equation}
v_{ebt}=\overline{v}_{eb\left(t+L\left(b\right)\right)}\text{ for each }e\in\mathcal{E},t\in\mathcal{T},b_{p}\in\mathcal{P}_{\text{stain}},
\end{equation}
where $\overline{v}_{eb\left(t+L\left(b\right)\right)}$ marks the end point of a batch and $b_{p}$ denotes the product contained in batch $b$. For a given edge $e$ and product $p_{0}\in\mathcal{P}_{\text{stain}}$ define the set of excluded batches as

\begin{equation}
E_{e,p_{0}}=\left\{ b\in e\text{ such that }b_{p}\in\mathcal{P}_{\text{stain}}\text{ and }b_{p}\neq p_{0}\right\}. 
\end{equation}
Then the exclusion condition prevents the placement of illegal sequences of batches reads
\begin{equation}
    \sum_{b\in E_{e,p_{0}}}v_{ebt}+\overline{v}_{eb_{0}t}\leq1\text{ for each }e\in\mathcal{E},t\in\mathcal{T},b_{0}\in e\text{ such that }\left(b_{0}\right)_{p}\in\mathcal{P}_{\text{stain}}.
\end{equation}
To enforce the flushing of the pipeline, we define
\begin{equation}
F_{e,b_{0}}=\left\{ b\in e\text{ such that }b_{p}\in\mathcal{P}_{\text{flush}}\text{, }b_{r}=\left(b_{0}\right)_{r}\text{ and }b_{V}\geq\left(b_{r}\right)^{V}\right\} , 
\end{equation}
where $b_{V}$ and $b_{r}$ are the batch size and regime of batch $b$, respectively. The flushing is then enforced by
\begin{equation}
    \overline{v}_{eb_{0}t}-v_{eb_{0}t}-\sum_{b\in F_{e,b_{0}}}v_{ebt}\leq0\text{ for each }e\in\mathcal{E},t\in\mathcal{T},b_{0}\in e\text{ such that }\left(b_{0}\right)_{p}\in\mathcal{P}_{\text{stain}}.
\end{equation}
To reduce the number of constraints, one can replace $e$ by $e_{\text{init}}$ since the placement of the batches in $e_{\text{tansit}}$ and $e_{\text{final}}$ is fully determined by $e_{\text{init}}$ (due to the routes constraints -- see Section \ref{sec_pregs_routes}).

\subsection{Regime Exclusions}\label{sec_regime_excl}
Technical restrictions dictate the need to exclude the simultaneous use of pumping regimes in certain scenarios.

To allow the simultaneous pumping of only one regime from a group of regimes at a given time, we define
\begin{equation}
R_{i}=\left\{ r_{1},r_{2},\ldots,r_{l_{i}}\right\} 
\end{equation}
where $i\in\left\{ 1,2,\ldots,k\right\}$. The sets $R_{i}$ represent groups of regimes on which we wish to impose mutual exclusions. Additionally, define
\begin{equation}
B_{R_{i}}=\left\{ \left(b,e\right):b\in e_{\text{init}}\text{ for some }e\in\mathcal{E},b_{r}\in R_{i}\right\}  
\end{equation}
as the set of all initial batches that relate to the pumping regimes in group $i$.
Then, the constraints which prevent the simultaneous placement of the batches that exclude each other read
\begin{equation}\label{regime_excl}
\sum_{\left(b,e\right)\in B_{R_{i}}}\sum_{t^{\prime}\in\left\{ t,t+1,\ldots,t+L\left(b\right)\right\} }v_{ebt^{\prime}}\leq1\text{ for each }t\in\mathcal{T}_{b},i\in\left\{ 1,2,\ldots,k\right\},
\end{equation} 
where $\mathcal{T}_{b}$ represents the time $\mathcal{T}$ horizon shortened in such a way that the right sum makes sense. Notice that the constraint \eqref{regime_excl} is formulated only for the initial batches as the exclusion behavior propagates due to the routes conditions detailed in Section \ref{sec_pregs_routes}.

\subsection{Storage Site  Capacity Constraints}\label{sec_store_s_const}
The Sections \ref{sec_bin_pack} - \ref{sec_regime_excl} described the dynamics of the transport through the network. The constraints introduced thus far do not track the storage capacities at each storage site, which can lead to the following unwanted behaviors:
\begin{itemize}
    \item a batch might be sent out from a site without it being physically present at the site,
    \item a batch might be received at a site without there being sufficient capacity available.
\end{itemize}
To remedy this, we introduce two occupancy variables:
\begin{itemize}
    \item $c_{s,p,t}^{\text{upper}}$ as the "blocked occupancy", which gets incremented by the given batch volume (given in discrete volume units) when a batch of product $p$ starts "pouring" into the tanks at site $s$ at time $t$.
    \item  $c_{s,p,t}^{\text{lower}}$ as the "on stock occupancy", which increases by the given batch volume (given in discrete volume units) when a batch of product $p$ is present at site $s$.
\end{itemize}
After these are computed, they will be bounded from above and below, respectively, which excludes the aforementioned "pathological" behaviors.
Note that the quantities $c_{s,p,t}^{\text{upper}}$ and $c_{s,p,t}^{\text{lower}}$ are given in discrete units, which can be used to compute the physical volume. Either way, we refer to the volume (or units of volume) of product for a given batch $b\in\mathcal{E}$ as $b_{V}$. Lastly, note that capacity tracking is only relevant for storage sites; therefore, we restrict our attention to $\mathcal{S}_{\text{store}}$ only.

To group the relevant incoming and outgoing batches for a given site define
\begin{equation}
B_{s,p}^{\text{in}}=\left\{ \left(b,e\right):e\in\mathcal{E},e_{\text{dest}}=s,b\in e_{\text{final}},b_{p}=p\right\}  
\end{equation}
\begin{equation}
B_{s,p}^{\text{out}}=\left\{ \left(b,e\right):e\in\mathcal{E},e_{\text{orig}}=s,b\in e_{\text{init}},b_{p}=p\right\} . 
\end{equation}
Additionally, the amount of product present at a site can change due to externalities (independent of the pipeline transport). To model this, we define $c_{s,p,t}^{\text{base}}$ as the (base) amount of product $p$ at site $s$ at time $t$.

With this in place, we compute $c_{s,p,t}^{\text{upper}}$ and $c_{s,p,t}^{\text{lower}}$ as follows
\begin{equation}\label{eq-upper-cap}
\sum_{\left(b,e\right)\in B_{s,p}^{\text{in}}}\sum_{t'\leq t}b^{V}v_{ebt'}-\sum_{\left(b,e\right)\in B_{s,p}^{\text{out}}}\sum_{t'\leq t-L\left(b\right)}b^{V}v_{ebt'}+c_{s,p,t}^{\text{base}}=c_{s,p,t}^{\text{upper}},\forall p\in\mathcal{P},\forall s\in\mathcal{S}_{\text{store}},\forall t\in\mathcal{T},
\end{equation}
\begin{equation}\label{eq-lower-cap}
\sum_{\left(b,e\right)\in B_{s,p}^{\text{in}}}\sum_{t'\leq t-L\left(b\right)}b^{V}v_{ebt'}-\sum_{\left(b,e\right)\in B_{s,p}^{\text{out}}}\sum_{t'\leq t}b^{V}v_{ebt'}+c_{s,p,t}^{\text{base}}=c_{s,p,t}^{\text{lower}},\forall p\in\mathcal{P},\forall s\in\mathcal{S}_{\text{store}},\forall t\in\mathcal{T},
\end{equation}

where the summations $\sum_{t'\leq t-L\left(b\right)}$ and $\sum_{t'\leq t-L\left(b\right)}$ are understood to be empty if $t-L\left(b\right)<0$. The quantity $c_{s,p,0}^{\text{base}} \geq 0$ represents the initial occupancy at a site and $c_{s,p,t}^{\text{base}}$ for $t>0$ is given based on the predicted outtake and may be negative. The maximal capacity per site $s$, product $p$ and time $t$ is denoted $c_{s,p,t}^{\text{max}}$. Then the constraints
\begin{equation}
c_{s,p,t}^{\text{upper}}\leq c_{s,p,t}^{\text{max}},\quad\forall p\in\mathcal{P},\forall s\in\mathcal{S}_{\text{store}},\forall t\in\mathcal{T},
\end{equation}
\begin{equation}
0\leq c_{s,p,t}^{\text{lower}},\quad\forall p\in\mathcal{P},\forall s\in\mathcal{S}_{\text{store}},\forall t\in\mathcal{T}\label{min_cap}
\end{equation}
ensure that we never overstock a storage site and never "ship out" a batch that is not present at a site. In case there are minimal requirements $c_{s,p,t}^{\text{min}}$ on the amount of product at a site, (\ref{min_cap}) may be replaced by
\begin{equation}
    c_{s,p,t}^{\text{min}}\leq c_{s,p,t}^{\text{lower}},\quad\forall p\in\mathcal{P},\forall s\in\mathcal{S}_{\text{store}},\forall t\in\mathcal{T}
\end{equation}
Note that the decision variables $c_{s,p,t}^{\text{upper}}$, $c_{s,p,t}^{\text{lower}}$ will also be used in Section \ref{sec_obj_function} to achieve a desired state at the end of the planning horizon (a desired distribution of products throughout the network).

Note that the values $c_{s,p,t}^{\text{min}}$ and $c_{s,p,t}^{\text{max}}$ are based on the operational limits for a given storage site. Due to the formulation of the update of storage capacities \eqref{eq-upper-cap} and \ref{eq-lower-cap} they must include an additional buffer that includes the counting error that is at most the volume of the pipes mediating the transport.

\subsection{Outages, Scheduling and Semi-product Transport}

\subsubsection{Outages}\label{ses-outages}
Sometimes, a pipeline, tank or other component of the network is out of order for a limited amount of time. This naturally influences the optimization process and requires care.
Since the maximum capacities per product at a site are handled by $c_{s,p,t}^{\text{max}}$, the tank outages are reduced to modify this variable for a period of time i.e.

\begin{equation}
    c_{s,p,t}^{\text{max}}=\hat{c}_{s,p,t}^{\text{max}}-c,\quad\forall t\in\left\{ t_{1},\ldots,t_{k}\right\}, 
\end{equation}
where $c$ is the reduction of capacity at a site due to the tank outage, $\left\{ t_{1},\ldots,t_{k}\right\}$ represents the duration over which the tank is out of order and $\hat{c}_{s,p,t}^{\text{max}}$ is the capacity before considering the outage.

Outages related to transport are defined by a set of batches and times over which the placement of batches is forbidden. More succinctly, let 
\begin{equation}
    B^{\text{exclusion}}=\left\{ \left(e_{1},b_{1}\right),\ldots,\left(e_{k},b_{k}\right)\right\} 
\end{equation}
be the set of relevant batches and $\left\{ t_{1},\ldots,t_{k}\right\} $ be the times over which the batch placement should not be possible. Then the associated constraints read
\begin{equation}
    v_{ebt}=0\text{ for each }\left(e_{1},b_{1}\right)\in B^{\text{exclusion}},t\in\left\{ t_{1},\ldots,t_{k}\right\}.
\end{equation}
Using the aforementioned constraint, it is possible to model regime outages, edge outages (when the pipeline needs to get cleaned, for instance), or product-specific outages (an inability to transport a particular product to a given site during a time interval). 

\subsubsection{Scheduling}\label{sec_scheduling}
There are situations in which only a particular amount of product is allowed to pass through an edge (or a group of edges) during a subset of the planning horizon $\mathcal{T}$. This occurs, for example, when considering the production limitations of a refinery $s$. In this case, all edges $e$ for which $e_{\text{orig}}=s$ will be subject to a cumulative product limit throughout the relevant time window. 

Let $E=\left\{ e_{1},\ldots,e_{k}\right\} $ be the set of edges on which we wish to impose the limit. Consider a time window $\hat{\mathcal{T}}\subset\mathcal{T}$. Let $L$ denote the limit for product $p$ (in units of volume) during $\hat{\mathcal{T}}$. Define
\begin{equation}
    P_{E}=\left\{ \left(e,b\right):e\in\mathcal{E},b_{r}=p\right\}, 
\end{equation}
then the constraint reads

\begin{equation}
\sum_{\left(e,b\right)\in P_{E}}\sum_{t\in\hat{\mathcal{T}}}b^{V}v_{ebt}\leq L.
\end{equation}

\subsubsection{Fixed-schedule Transport}
On occassion, it is necessary to carry out the transport of products between refineries at fixed times. These products never enter any storage site within the network and are therefore handeled as an outage of fixed length.

\subsection{Objective Function}\label{sec_obj_function}
The following considerations motivate the formulation of the cost function:
\begin{itemize}
    \item During the time horizon, the amount of product coming for refineries is maximized up to a certain quota.
    \item Pumping cost should be minimized.
    \item (Optional) The distribution of products at the end of the planning horizon throughout the network should tend towards a defined "ideal state".
    \item (Optional) Sometimes, it is necessary to base a plan on an existing one. In this case, a term is added that attracts the state to the former plan.
\end{itemize}
Before the components of the objective function are discussed, a constraint which enforces a maximal quota on the intake from refineries is formulated. These constraints complete the bin-packing problem (\ref{box_spatial_count_1})-(\ref{box_overlap_exclusion}) is added. The consequence of setting this kind of quota on the intake from refineries allows a hierarchical setup for the cost function in which the intake from refineries can be prioritized over the optimization of pumping costs (see Section \ref{sec_num_study}).

\subsubsection{Nominations and Related Additional Constraints}
Let the maximum amount of product intake for a given refinery $s\in\mathcal{S}_{\text{ref}}$ defined by
\begin{equation}
    N_{s}=\left(\left(p_{1},V_{1}\right),\left(p_{2},V_{2}\right),\ldots,\left(p_{k},V_{k}\right)\right),\label{nomination_form}
\end{equation}
where $s$ is a refinery and the list of ordered pairs $\left(p_{1},V_{1}\right),\left(p_{2},V_{2}\right),\ldots,\left(p_{k},V_{k}\right)$ gives the information on each product and maximal amount.  Let the maximum product intake of the product $p$ from the refinery $s$ subject to nomination $N_{s}$ be denoted $V_{s,p}$

For a given refinery $s$, define the set of edges connected to the refinery as
\begin{equation}
    \mathcal{E}_{s}=\left\{ e\in\mathcal{E}:e_{\text{orig}}=s\right\} \label{E_N}
\end{equation}
Furthermore, define the set of all boxes on edge $e$ corresponding to product $p_{0}$ as
\begin{equation}
    e_{p_{0}}=\left\{ \left(e,b\right):b\in e,b_{p}=p_{0}\right\}. \label{e_p0}
\end{equation}
The additional constraint relevant to the nomination $N_{s}$ then reads
\begin{equation}\label{nomination-eq}
    \sum_{e\in \mathcal{E}_{s}}\sum_{\left(e,b\right)\in e_{p_{i}}}\sum_{t\in\mathcal{T}}b^{V}v_{ebt}\leq V_{s,p}\text{ for each }p\in N_{s},
\end{equation}
where a slight abuse of notation was committed. If multiple refineries supply the product to the network, the constraint \eqref{nomination-eq} is applied for each of these.

\subsubsection{Structure of the Objective function}\label{sec_structure_of_obj}
As mentioned previously, the objective function to by maximized has multiple components and reads
\begin{equation}
J=\alpha J_{N}+\beta J_{D}+\gamma J_{P}+\theta J_{C}\label{cost_function}
\end{equation}
where $J_{N},J_{D},J_{P} \text{ and } J_{C}$ are the refinery inflow, distributional, re-planning and pumping cost components, respectively. The coefficients $\alpha,\beta,\gamma,\theta\geq0$ are used to establish the relative importance of each of the components of the objective function. The following sections detail each of the components of (\ref{cost_function}).

\subsubsection{The Refinery Inflow Component of the Objective Function}
The component $J_{N}$ evaluates the utility of fulfilling a given maximal inflow from the refineries. To show how it is constructed, assume that $N_{s}$ is of the form (\ref{nomination_form}). Using the definitions of $E_{n}$ and $e_{p_{0}}$ given by (\ref{E_N}) and (\ref{e_p0}) we define
\begin{equation}
    J_{N_{s}}=\sum_{p_{i}\in N_{s}}\sum_{e\in \mathcal{E}_{s}}\sum_{\left(e,b\right)\in e_{p_{i}}}\sum_{t\in\mathcal{T}}\eta_{p_{i}}b^{V}v_{ebt},
\end{equation}
where $\eta_{p_{i}}$ attributes relative importance to each of the products in the nomination. Variables $\eta_{p_{i}}$ and the regularization coefficient $\alpha$ can be set so that $\alpha\eta_{p_{i}}$ represents the monetary value of a unit of product $p_{i}$. Using this scaling, it is possible to directly fine-tune this against the pumping regime costs expressed by $J_{C}$ (Discussed in Section \ref{sec_nom_comp}). Alternatively, one can set the coefficients $\alpha$ and $\gamma$ so that one has priority over the other.

\subsubsection{The Distributional Component of the Objective Function}

To ensure that market demands are met, it is meaningful to optimize in an extended time horizon. If this is not computationally tractable, one can resort to prompting a certain end state in the cost function. Two methods of doing this are explained in the following.

In the first method, an exact occupancy is defined for each site at final time. More concisely, let $c_{s,p}^{\text{opt}}\geq 0$ and $k_{s,p}>0$ denote the optimal occupancy and weight coefficient of the product $p$ at site $s$, respectively. The distributional component of the cost function then reads

\begin{equation}
J_{D}=-\sum_{s\in\mathcal{S}}\sum_{p\in\mathcal{P}}k_{s,p}\left|c_{s,p,\mathcal{T}_{\text{final}}}^{\text{lower}}-c_{s,p}^{\text{opt}}\right|
\end{equation}

In the second method, no target occupancy is used, instead, the following variables are defined
\begin{equation}
    k_{s,p}^{\prime}=\begin{cases}
k_{s,p} & \text{if it is desired to have as much product p at site s as possible}\\
-k_{s,p} & \text{if it is desired to have as little product p at site s as possible}.
\end{cases}
\end{equation}
Using this definition, the alternative distributional component reads
\begin{equation}
    J_{D}^{\prime}=\sum_{s\in\mathcal{S}}\sum_{p\in\mathcal{P}}k_{s,p}^{\prime}c_{s,p,\mathcal{T}_{\text{final}}}^{\text{lower}}.
\end{equation}
Note that either of these methods can be applied to subsets of sites and combined with each other.

\subsubsection{The Planning Component of the Objective Function}
If circumstances change, one may want to re-plan within the time horizon $\mathcal{T}$, using a known schedule. This should lead to a new plan that also takes into account the previous plan. Let $v_{ebt}^{\text{prev}}$ denote the decision variables associated with the previous plan and let 
\begin{equation}
    \vartheta^{\text{prev}}=\left\{ \left(e,b,t\right):v_{ebt}^{\text{prev}}\neq0\right\}, 
\end{equation}
then the planning component of the objective function reads
\begin{equation}
J_{P}=-\sum_{\left(e,b,t\right)\in\vartheta^{\text{prev}}}\left|v_{ebt}^{\text{prev}}-v_{ebt}\right|^{2}.\label{planning_comp}
\end{equation}
It may also be desirable to only define the planning component on a subset of the original plan, in this case, just replace $\vartheta^{\text{prev}}$ in (\ref{planning_comp}) by a subset of the plan.

Occasionally, a portion of the plan might have already been executed. In this case we define $\vartheta^{\text{executed}}\subset\vartheta^{\text{prev}}$ as the subset of batches that have already been executed and add the constraints
\begin{equation}
    v_{ebt}=1\text{ for each }\left(e,b,t\right)\in\vartheta^{\text{executed}}.
\end{equation}
\subsubsection{The Pumping Regime Cost Component of the Objective Function}\label{sec_nom_comp}
A cost is associated with any batch scheduled. This is captured by the regime cost component of the objective function. Let $b_{\text{cost}}$ be the cost associated with the scheduling of batch $b$. Then the regime cost component reads

\begin{equation}
    J_{R}=-\sum_{e\in\mathcal{E}}\sum_{b\in e}\sum_{t\in\mathcal{T}}b_{\text{cost}}v_{ebt}.
\end{equation}

\section{Notes on the Algorithmic Implementation}\label{sec_const_handl}
When formulating the model, tests were performed using the SCIP optimization suite. These tests included path graph and tree graph topologies and led to choosing particular formulations for given constraints (see Section \ref{sec_MILP_form}). These experiments also provided insight into the scaling of the formulation when using SCIP. 

By far the greatest increase in computational cost was observed when adding the tank capacity constraints detailed in Section \ref{sec_store_s_const}. This is hardly a surprise, as they interconnect all of the bin-packing problems, which were only sparsely connected before these constraints were imposed (see \eqref{eq-lower-cap}, \eqref{eq-upper-cap}). Furthermore, during the same set of scaling tests, it was observed that many of the constraints of the type \eqref{eq-lower-cap}, \eqref{eq-upper-cap} were not active for a considerable amount of the test instances. For this reason, we make use of the constraint handler feature of the SCIP optimization suite, which allows us to add constraints dynamically. The effects of adding the constraints in this fashion are studied in Section \ref{sec_scaling_scip}

\section{A Numerical Study}\label{sec_num_study}

In this section, several scaling tests are presented. Due to the large amount of recent publications, which consider a path graph topology, all of the experiments are performed on a path graph. This allows for a more direct comparison with other formulations, which is provided in Section \ref{sec_scaling_comparison}. The formulation has also been applied to other topologies and is currently used to provide monthly planning for the Czech national pipeline operator. Reporting on these will require a more comprehensive numerical study, which we plan to provide in a future publication.

Some of the experiments are performed using our custom algorithmic solution implemented in SCIP, which makes use of the constraint handler function (see Section \ref{sec_const_handl}), and vanilla SCIP. The aforementioned tools are open-source, but typically exhibit lesser performance compared to commercial solvers. Due to this performance gap, Gurobi Optimizer is used to benchmark the overall performance of the formulation. In all the numerical experiments, a two-product scenario is considered, in which one of the products is a staining product and the other flushing and the time horizon is discretized in increments on one hour, which is a standard setting in many benchmarks \cite{Li2021}. Note that the scaling of the physical quantities does not necessarily correspond to any real-world scenario.

\subsection{General Setup}

All of the following simulations take place on a path graph of length $l$. Each path graph has a single refinery at the very end of the graph and $l-1$ storage nodes that are represented by the other vertices of the graph.
All numerical experiments are performed on an Apple Macbook Pro laptop equipped with 16 GB of RAM and an Apple M2 Pro chip.

\subsection{Scaling using SCIP}\label{sec_scaling_scip}

In this section, the performance of the formulation using the open source solver SCIP is studied. Note that the size of instances is comparable to the ones studied using commercial solvers for other formulations \cite{Li2021, TSUNODAMEIRA2021105143}. Section \ref{sec_scaling_comparison} contains a more leveled comparison between our formulation and the state of the art.

Table \ref{tab_settings_of_dem_path} summarizes the settings which apply to all the experiments. Notably, the initial occupancies of each of the sites increase linearly with the distance from the refinery, and the batch sizes are set to values close to what may occur in practice. The pump speed is set so that precisely integer batch sizes (in time units) are obtained. All of the experiments focus on scaling with respect to the length of the path graph, which allows a head-to-head comparison with other formulations.

Outside of the settings summarized in Table \ref{tab_settings_of_dem_path} the experiments assume the additional settings listed in Table \ref{tab_additional_settings_path}. Lastly, different cost functions settings are considered. A summary of these can be reviewed in Table \ref{tab_different_cost_function_modes}. Lastly, we analyze solution times across three different solver setups: vanilla SCIP, SCIP with constraint handler (detailed in Section \ref{sec_const_handl}) and Gurobi Optimizer.

The first of the three experiments focuses on the scaling with respect to supply and demand requirements (SD). In this set of instances, the outtakes at each site are enforced by constraints, while the refinery supply is maximized by means of the cost function (see Table \ref{tab_different_cost_function_modes}). In this case, no restrictions on the solution time are given and the optimum has been reached in all cases, as documented by Table \ref{tab_scaling_experiments_1}. Inspecting the figures shows the following two features:
\begin{itemize}
    \item The constraint handler implementation may speed up or slow down the process based largely on the amount of callbacks to the constraint handler.
    \item The scaling figures are notably better than those reported in the recent literature. We refer to Section \ref{sec_scaling_comparison} for full discussion.
\end{itemize}

The next set of experiments, whose details are laid out in Table \ref{tab_scaling_experiments_2} is concerned with optimizing the costs as well as meeting the supply and demand SDC (see Table \ref{tab_different_cost_function_modes}). In contrast with the previous case, the optimal solution is not set as a stop condition. Instead, the computation is stopped when either the optimality gap falls below $10^{-3}$ or a computational time of $30$ minutes is reached. Note that because of the structure of the cost function, the refinery supply is always prioritized over the minimization of pumping (energy savings), which leads to the available supply being extracted in all the cases. In all of these instances, we report basic SCIP having better scaling properties than the constraint handler implementation due to the large number of callbacks invoked when solving the problems at hand. Nevertheless, the problem can be solved to a satisfactory error for up to 7 sites in 30 minutes. Due to the proportion coefficients of the cost function, we may directly compare the experiments of Table \ref{tab_scaling_experiments_1} with the present set of experiments to help us gauge the effect of the multi-objective optimization and the duality gap. Table \ref{tab_different_cost_function_modes} shows the pumping cost reductions that may be achieved by considering the SDC cost function setting for settings B. 

So far, all of the scaling tests, have made use of SCIP, which is an open source solver. Most of the recent publications make use of commercial solvers such as IBM ILOG CPLEX or Gurobi Optimizer. To provide more accurate information on the scalability of our formulation especially for larger instances, we provide tests that make use of Gurobi Optimizer and solve large instances to complete the picture, leading to the final scaling comparison of Section \ref{sec_scaling_comparison}.

\begin{table}
  \centering
  \caption{Cost function settings for each of the experiments}
  \label{tab_different_cost_function_modes}
\begin{tabular}{p{5cm}|p{5cm}|p{5cm}}
 \toprule
 Experiment label& Description & The cost function setting according to \eqref{cost_function}  \\
 \midrule
Supply demand only (SD) & Leaving only the refinery component non-zero satisfies the outtakes at the sites while maximizing the intake from the refinery. & $\alpha\neq0, \theta=0,\beta=0 \text{ and }\gamma=0$ \\
Supply demand and pumping regime costs (SDC) & Modifying the setup of SD by changing the pumping regime costs coefficient to a non-zero value allows for the intakes from the refinery to be prioritized, while finding a more efficient way to pump. & $\alpha=5, \theta=3*10^{-3},\beta=0 \text{ and }\gamma=0$ \\
\bottomrule
\end{tabular} 
\end{table}

\begin{table}[H]
  \centering
  \caption{Additional settings applied to the path graph experiments}
  \label{tab_additional_settings_path}
\begin{tabular}{ p{5cm}|p{4cm}|p{4cm}  }
 \toprule
 Experiment label& Refinery intake per product in batches &Time horizon \\
 \midrule
A - small refinery intake & 10/10  &480 time units (20 days)\\
B - moderate refinery intake &   15/15  & 576 time units (24 days)\\
C - large refinery intake & 20/20 & 576 time units (24 days) \\
\bottomrule
\end{tabular} 
\end{table}

\begin{table}[H]
  \centering
  \caption{Settings common to all of the experiments}
  \label{tab_settings_of_dem_path}
\begin{tabular}{p{6cm}|p{4cm}|p{3cm}}
 \toprule
 Parameter& Physical value &Value in units \\
 \midrule 
 Initial storage occupancy of the flushing / staining product at the first site   & 6976.8 / 7792.8 $m^{3}$  &120/120\\
 Storage occupancy increase per site for flushing / staining products &   581.4 / 649.4  $m^{3}$ & 10 / 10 \\
 Storage site outtake per product & 581.4 / 649.4  $m^{3}$ & 10 / 10 \\
 Amount of product per batch & 5814 / 2857.36 $m^{3}$ &  100 / 44 \\
 Pumping speed per product & 969 / 952.45 $\frac{m^{3}}{t}$ & - \\
 Time to pump a batch into the pipeline & 6 / 3 hours & 6 / 3 time units \\
\bottomrule
\end{tabular} 
\end{table}

\begin{table}[H]
  \centering
  \caption{Scaling experiments of the supply demand experiment (SD). The top and bottom tables being the SCIP constraint handler and vanilla SCIP solution times, respectively. The optimum was achived in each of the cases. The solving times are given in seconds.}
  \label{tab_scaling_experiments_1}
\begin{tabular}{p{3cm}|p{2cm}|p{2cm}|p{2cm}}
 \toprule
 Number of vertices & setting A & setting B & setting C \\
 \midrule
 4 & 27.84 & 36.82 & 1221 \\
 5 & 106.9 & 249 & 3520.1 \\
 6 & 127.9 & 1882.9 & 3071.6 \\
 7 & 216.8 & 521.4 & 6929.9 \\
 8 & 4698.1 & infeasible & 378.3 \\
\midrule
 4 & 211.4 & 360.9 & 331.3 \\
 5 & 364.8 & 402.5 & 462.8 \\
 6 & 154.1 & 548.8 & 533.9 \\
 7 & 446.9 & 762.2 & 567.1 \\
 8 & 466.9 & infeasible & 721.9 \\
\bottomrule
\end{tabular} 
\end{table}

\begin{table}[H]
  \centering
  \caption{Scaling experiments of the supply demand and pumping regimes costs experiment (SDC). The top and bottom tables being the constraint handler SCIP and vanilla SCIP solution times, respectively. The solving times are given in seconds and the optimality gap is noted in brackets. The computation is stopped if either a 30 minute duration is reached or the duality gap is below e-3.}
  \label{tab_scaling_experiments_2}
\begin{tabular}{ p{3cm}|p{3cm}|p{3cm}|p{3cm}  }
 \toprule
 Number of vertices & setting A & setting B & setting C \\
 \midrule
 4 & 75.15 (9.858e-4)& 475.6 (4.167e-4) & 1809.6 (0.008) \\
 5 & 137.1 (6.527e-4)& 1814.1 (0.005) & 1813.3 (0.009) \\
 6 & 276.7 (9.883e-4)& 1817.8 (0.003) & 1817.3 (0.007) \\
 7 & 459.8 (3.511e-5) & 1821.2 (0.001) & 1820.4 (inf) \\
 8 & 1817.4 (0.006) & infeasible & 1824.4 (inf) \\
\midrule 
 4 & 232.3 (3.454e-5) & 387.5 (6.407e-4) & 377.5 (7.301e-4)\\
 5 & 255.5 (3.477e-5) & 368.2 (6.436e-4) & 199.1 (4.649e-4)\\
 6 & 1812.2 (0.006) & 1817.5 (0.002) & 1821.5 (0.002)\\
 7 & 421.2 (5.620e-4) & 1821.2 (0.001) & 1822.4  (0.002)\\
 8 & 1817.5 (0.006) & infeasible & 1825 (inf) \\
\bottomrule
\end{tabular} 
\end{table}

\begin{table}[H]
  \centering
  \caption{A comparison of pumping costs when considering the SD and SDC cost function settings (see Table \ref{tab_different_cost_function_modes}). The comparison is made for setting B of Table \ref{tab_different_cost_function_modes} and the activity of a pumping regime is awarded with a cost of 1 per hour.}
  \label{tab_gap_comparison}
\begin{tabular}{p{3cm}|p{3cm}|p{3cm}|p{4cm}}
 \toprule
 Number of vertices & SDC & SD & improvement in percent \\
 \midrule
 4 & 1452& 1752 & 21 \\
 5 & 1608& 2340 & 46 \\
 6 & 1836& 2988 & 62 \\
 7 & 2064 & 3072 & 49 \\
\bottomrule
\end{tabular} 
\end{table}

\subsection{A Comparison of Scaling Properties of Multiple MILP Formulations}\label{sec_scaling_comparison}
As mentioned before, the results mentioned in Section \ref{sec_scaling_comparison} do not fully demonstrate the capabilities of the formulation. To arrive at a fair comparison, we solve larger instances using Gurobi Optimizer and compare the results with the formulations presented earlier, as surveyed by Li et al. \cite{Li2021}.

In the literature \cite{Liang16, Ghaithan20}, multi-product pipeline planning and scheduling problems are traditionally divided into short-term (4-14 days) and long-term (over 2 weeks), where the time horizon is divided into hour-long increments for the purposes of the MILP formulation. For the longer-term planning problem, there are formulations that yield computational results for up to 6 sites \cite{meira20, Li2021}. Others have studied situations, where they scale the product count upward, but end up providing results for 4 or 5 sites only \cite{Li2021, Liang16, meira20}, typically on shorter time horizons. The time limit is not universally agreed on, but generally falls in the range of 1--2 hours of computational time \cite{Li2021}.

To set the stage for this key scaling experiment, we consider the settings of Table \ref{tab_settings_of_dem_path} and a nomination, which is of 40/40 batches (see Table \ref{tab_additional_settings_path} for a comparison). Lastly, considering the SDC setting for the cost function (see Table \ref{tab_different_cost_function_modes}) we optimize for a very long horizon (up to two months) and a large network (up to 11 sites). The results of the experiments are detailed in Table \ref{tab_scaling_experiments_3}. It is immediately apparent that the present formulation exhibits favorable scalability. This is further confirmed by comparing the results of Table \ref{tab_scaling_experiments_3} to the current state of the art. The results of this comparison are presented in Table \ref{tab_big_scaling_comp}.

\begin{table}[H]
  \centering
  \caption{The scaling experiments for larger networks and longer time horizons using Gurobi Optimizer. The stop condition for each of the simulations if the optimality gap reaching below $10^{-4}$. The refinery intake for the instances marked by a * is doubled to 80/80 units for these instances to be feasible.}
  \label{tab_scaling_experiments_3}
\begin{tabular}{p{3cm}|p{3cm}|p{3cm} }
 \toprule
 Number of vertices & 31 days (744 time units) & 62 days (1488 time units). \\
 \midrule
 4 & -  & 144 \\
 5 & - & 1116* \\
 6 & -  & 452* \\
 7 & 858  & 906* \\
 8 & 414 & - \\
9 & 1695 & -\\
10 & 610 & -\\
11 & 692 & -\\
12 & 989 & -\\
\bottomrule
\end{tabular} 
\end{table}

\begin{table}[H]
  \centering
  \caption{Comparison with state of the art results present in literature. All of the publications below deal with a path graph topology (single line) network with a single source and outtakes at each of the storage nodes.}
  \label{tab_big_scaling_comp}
\begin{tabular}{ p{4cm}|p{2.5cm}|p{2.5cm}|p{2.5cm}  }
 \toprule
 Authors & Number of nodes & Time horizon in days & Solving time in seconds \\
 \midrule
Meira et al. \cite{meira20} & 6 & 30 & 1476 \\
Rejowski et al. \cite{pinto04} & 5 & 3.2 & 2800 \\
Xu et al. \cite{XU2021341} & 5 & 15 & 186 \\
Mostafaei et al. \cite{most20} & 4 & 8.3 & 2480 \\
\midrule
Proposed  formulation {scaled~w.r.t.~the~topology} & 12 & 31 & 989 \\
Proposed formulation {scaled~in~time} & 7 & 62 & 906 \\
\bottomrule
\end{tabular} 
\end{table}

\section{Conclusion}
A novel MILP formulation for scheduling commodity transport in a multi-product pipeline, based on the spatially-indexed bin packing problem \cite{Allen_2012}, was proposed. The favorable scaling properties of the aforementioned formulation are maintained even after all the operations-specific constraints are added (Sections \ref{sec_pregs_routes} - \ref{sec_scheduling}). Using a set of path-graph test instances, the scaling properties are compared to existing formulations and shown to be a remarkable improvement. More specifically, it is possible to solve problem instances with time horizons or sizes that extend beyond the current state of the art (see Table \ref{tab_big_scaling_comp}). Additionally, validating the proposed formulation is the recent deployment of an information system at the Czech national operator, which we plan to detail in a follow-up paper. 

\bibliographystyle{ieeetr}
\bibliography{main}

\end{document}